\documentclass{article}
\usepackage[latin9]{inputenc}
\usepackage{mathtools}
\usepackage{amsmath}
\usepackage{amsthm}
\usepackage{amssymb}

\makeatletter
\numberwithin{equation}{section}
  \theoremstyle{remark}
  \newtheorem*{rem*}{\protect\remarkname}
\theoremstyle{plain}
\newtheorem{thm}{\protect\theoremname}
  \theoremstyle{plain}
  \newtheorem{cor}{\protect\corollaryname}
 \theoremstyle{definition}
 \newtheorem*{defn*}{\protect\definitionname}
  \theoremstyle{plain}
  \newtheorem{prop}{\protect\propositionname}

\usepackage{amsfonts}
\usepackage{graphicx}
\providecommand{\U}[1]{\protect\rule{.1in}{.1in}}

\makeatother

  \providecommand{\definitionname}{Definition}
  \providecommand{\propositionname}{Proposition}
  \providecommand{\remarkname}{Remark}
\providecommand{\corollaryname}{Corollary}
\providecommand{\theoremname}{Theorem}

\begin{document}

\title{Existence of isotropic complete solutions of the \\ $\Pi$-Hamilton-Jacobi
equation}

\author{Sergio Grillo\\
{\small{}Instituto Balseiro, Universidad Nacional de Cuyo and CONICET}\\[-5pt]
{\small{}Av. Bustillo 9500, San Carlos de Bariloche}\\[-5pt] {\small{}R8402AGP,
República Argentina}}
\maketitle
\begin{abstract}
Consider a symplectic manifold $M$, a Hamiltonian vector field $X$
and a fibration $\Pi:M\rightarrow N$. Related to these data we have
a generalized version of the (time-independent) Hamilton-Jacobi equation:
the $\Pi$-HJE for $X$, whose unknown is a section $\sigma:N\rightarrow M$
of $\Pi$. The standard HJE is obtained when the phase space $M$
is a cotangent bundle $T^{*}Q$ (with its canonical symplectic form),
$\Pi$ is the canonical projection $\pi_{Q}:T^{*}Q\rightarrow Q$
and the unknown is a closed $1$-form $\mathsf{d}W:Q\rightarrow T^{*}Q$.
The function $W$ is called Hamilton's characteristic function. Coming
back to the generalized version, among the solutions of the $\Pi$-HJE,
a central role is played by the so-called \textit{isotropic complete
solution}s. This is because, if a solution of this kind is known for
a given Hamiltonian system, then such a system can be integrated up
to quadratures. The purpose of the present paper is to prove that,
under mild conditions, an isotropic complete solution exists around
almost every point of $M$. Restricted to the standard case, this
gives rise to an alternative proof for the local existence of a \textit{complete
family} of Hamilton's characteristic functions.
\end{abstract}

\section{Introduction}

In the recent few years, several extensions of the Hamilton-Jacobi
Theory have been developed. See for instance \cite{bmmp,car,pepin-noholo,lmm,hjp}.
In Ref. \cite{gp}, it was presented an extension to general dynamical
systems (on fibered phase spaces), which contains as particular cases
the previuos ones. Let us briefly review it in the restricted context
of Hamiltonian systems. Consider a symplectic manifold $\left(M,\omega\right)$
of dimension $d=2s$, a second manifold $N$ of dimension $k$ and
a surjective submersion $\Pi:M\rightarrow N$. Consider also a Hamiltonian
system with phase space $M$ and a Hamiltonian function $H$. According
to Ref. \cite{gp}, a \textbf{(global) complete solution} of the so-called
$\Pi$-\textbf{HJE} ($\Pi$-Hamilton-Jacobi Equation), for the Hamiltonian
vector field $X_{H}$, is a surjective local diffeomorphism $\Sigma:N\times\Lambda\rightarrow M$
such that 
\begin{equation}
\mathfrak{i}_{X_{H}^{\Sigma}}\Sigma^{*}\omega=\Sigma^{*}\mathsf{d}H\;\;\;\;\textrm{and}\;\;\;\Pi\circ\Sigma=p_{N},\label{hje}
\end{equation}
where $\Lambda$ is a third manifold of dimension $d-k$, the vector
field $X_{H}^{\Sigma}\in\mathfrak{X}\left(N\times\Lambda\right)$
is given by
\[
X_{H}^{\Sigma}\left(n,\lambda\right)=\left(\Pi_{*}\left(X_{H}\left(\Sigma\left(n,\lambda\right)\right)\right),0\right),
\]
and $p_{N}:N\times\Lambda\rightarrow N$ is the canonical projection
onto the first factor. Naturally, a \textbf{local complete solution}
of the $\Pi$-HJE for $X_{H}$, \textbf{along an open subset} $U$,
is a complete solution of the $\left.\Pi\right|_{U}$-HJE for $\left.X_{H}\right|_{U}$.
(Here, we are seeing $\left.\Pi\right|_{U}$ as a fibration onto $\Pi\left(U\right)$
and $\left.X_{H}\right|_{U}$ as a vector field on $U$). For each
$\lambda\in\Lambda$, the function 
\[
\sigma_{\lambda}:N\rightarrow M:n\mapsto\Sigma\left(n,\lambda\right)
\]
is a section of $\Pi$ and it is called a \textbf{partial solution}
(or simply a solution) of the $\Pi$-HJE. On the other hand, $\Sigma$
is said to be \textbf{isotropic} if 
\[
\sigma_{\lambda}^{*}\omega=0,\;\;\;\forall\lambda\in\Lambda,
\]
or, equivalently, if each vector subspace 
\begin{equation}
\Sigma_{*}\left(T_{n}N\times\left\{ 0\right\} _{\lambda}\right)\subseteq T_{\Sigma\left(n,\lambda\right)}M\label{ie}
\end{equation}
is isotropic with respect to $\omega$. Note that, in such a case,
we must have $s\geq k$. When $k=s$, above linear spaces are Lagrangian,
and consequently $\Sigma$ is said to be \textbf{Lagrangian}. 

The standard Hamilton-Jacobi Theory (see for instance \cite{am} and
\cite{car}) corresponds to the case in which:
\begin{itemize}
\item $M$ is a cotangent bundle $T^{*}Q$, 
\item $\omega$ is the canonical symplectic form $\omega_{Q}$ of $T^{*}Q$, 
\item $\Pi$ is the canonical projection $\pi_{Q}:T^{*}Q\rightarrow Q$,
\end{itemize}
and its complete solutions are precisely the Lagrangian complete solutions
of the $\pi_{Q}$-HJE. More precisely (see Ref. \cite{gp}), in this
case the complete solutions are locally given by the formula $\Sigma\left(q,\lambda\right)=\mathsf{d}W_{\lambda}\left(q\right)$,
where each function $W_{\lambda}$, called \textit{Hamilton's characteristic
function}, satisfies the \textit{standard} (time-independent) Hamilton-Jacobi
equation:
\begin{equation}
\mathsf{d}\left(H\circ\mathsf{d}W_{\lambda}\right)=0.\label{ehje}
\end{equation}
We shall say that the functions $W_{\lambda}$'s give a \textit{complete
family of }Hamilton's characteristic functions.

The importance of the isotropic complete solutions, as proven in \cite{g,gp}
(generalizing a well-known result of the standard theory), lies in
that fact that if we know one of them for a given Hamiltonian system,
then such a system is exactly solvable. Actually, in order to ensure
exact solvability, it is enough to know a local solution (instead
of a global one) around every point of the phase space. 

It was shown in \cite{gp} (see Theorem 4.15 there) that a local complete
solution of the $\Pi$-HJE always exists around every point $m\in M$
such that $X_{H}\left(m\right)\notin\mathsf{Ker}\Pi_{*,m}$. However,
nothing has been said about the existence of isotropic complete solutions.
In this paper we fill in this gap, showing that a local isotropic
complete solution does exist around every point $m\in M$ such that,
besides the condition $X_{H}\left(m\right)\notin\mathsf{Ker}\Pi_{*,m}$,
the subspace $\mathsf{Ker}\Pi_{*,m}\subseteq T_{m}M$ is co-isotropic
with respect to $\omega$. Moreover, if no fibration is previously
fixed, we show that around every point $m$ such that $X_{H}\left(m\right)\neq0$
there always exist a fibration $\pi$ and a local isotropic complete
solution of the $\pi$-HJE for $X_{H}$.

It is worth mentioning that our results do not mean that every Hamiltonian
system is exactly solvable around such points, because in order to
solve its equations of motion is not enough to ensure the existence
of an isotropic complete solution, but we need to have a concrete
expression of one of it. Thus, such results are mainly of a theoretical
nature, and can be interpreted as follows: every Hamiltonian system
is ``potentially'' exactly solvable\textit{ }with the aid of an
isotropic complete solution (following the procedure described in
Refs. \cite{g,gp}).

\bigskip{}

The organization of the paper is as follows. In Section 2 we recall
the duality between complete solutions and first integrals, because
we shall use it to prove our main results. In Section 3.1 we prove
the annunciated existence theorem and, as an immediate corollary,
we give in Section 3.2 a novel proof for the local existence of a
complete family of Hamilton's characteristic functions. As another
corollary, and only for completeness, we show in Section 3.3 that
every Hamiltonian system has a set of first integrals as those appearing
in the definition of a \textbf{non-commutative} or Mischenko-Fomenko
\textbf{integrable system} \cite{mf} (see also \cite{j} for a review
on the subject), defined around every non-critical point of the Hamiltonian
function. In other words, we prove, as a corollary of our main result,
that every Hamiltonian system is ``potentially'' non-commutative
integrable around almost every point of its phase space. (It is worth
mentioning that this fact can also be proven as a direct consequence
of the Carathéodory-Jacobi-Lie theorem \cite{lm}). Finally, in Section
3.4, using again the above mentioned duality, we show that around
every non-critical point of a Hamiltonian system there always exist
a fibration and a local isotropic complete solution related to it.

\bigskip{}

We assume that the reader is familiar with the basic concepts of Differential
Geometry (see \cite{boo,kn}), and with the basic ideas related to
Hamiltonian systems in the context of Symplectic Geometry (see \cite{am,lm,mr}).
We shall work in the smooth (i.e. $C^{\infty}$) category, focusing
exclusively on finite-dimensional smooth manifolds.

\section{The \textit{complete solutions - first integrals} duality}

\label{dual}

To show our main result, we shall use the duality between (isotropic)
complete solutions and (isotropic) first integrals stablished in \cite{gp}.
Let us recall such a duality for the case of Hamiltonian systems.
Following the same notation as in the Introduction, consider a symplectic
manifold $\left(M,\omega\right)$ of dimension $d=2s$, a function
$H:M\rightarrow\mathbb{R}$ and a fibration $\Pi:M\rightarrow N$
(i.e. a surjective submersion), with $\mathsf{dim}N=k$.

We shall say that a submersion $F:M\rightarrow\Lambda$ is a \textbf{first
integrals submersion} if 
\begin{equation}
\mathsf{Im}X_{H}\subseteq\mathsf{Ker}F_{*}.\label{fi}
\end{equation}
We adopt the following convention for defining $X_{H}$:
\begin{equation}
\mathsf{i}_{X_{H}}\omega=\mathsf{d}H.\label{conv}
\end{equation}

\begin{rem*}
Note that, if $\Lambda=\mathbb{R}^{l}$, the components $f_{1},...,f_{l}:M\rightarrow\mathbb{R}$
of $F$ define a set of $l$ (functionally) independent first integrals,
in the usual sense, for the Hamiltonian system defined by $H$. 
\end{rem*}
Also, we shall say that $F$ is \textbf{transverse to} $\Pi$ if 
\begin{equation}
TM=\mathsf{Ker}\Pi_{*}\oplus\mathsf{Ker}F_{*}.\label{tr}
\end{equation}
Finally, the map $F$ is said to be \textbf{isotropic} if 
\begin{equation}
\mathsf{Ker}F_{*}\subseteq\left(\mathsf{Ker}F_{*}\right)^{\bot},\label{is}
\end{equation}
where ``$\bot$'' indicates the symplectic orthogonal w.r.t. $\omega$.
Of course, in such a case $l\geq s$. When $l=s$, above inclusion
reduces to an equality and the map $F$ is say to be \textbf{Lagrangian}.

For later convenience, we shall say that an isotropic submersion $F$
is \textbf{symplectically complete} if $\left(\mathsf{Ker}F_{*}\right)^{\bot}$
is an integrable distribution.

\bigskip{}

It was shown in \cite{gp} that, given an isotropic complete solution
$\Sigma:N\times\Lambda\rightarrow M$ of the $\Pi$-HJE for $H$ {[}see
the Eqs. \eqref{hje} and \eqref{ie}{]}, we can construct around
every point of $M$ a neighborhood $U$ and a submersion $F:U\rightarrow\Lambda$
such that
\begin{itemize}
\item $\mathsf{Im}\left.X_{H}\right|_{U}\subset\mathsf{Ker}F_{*}$ (first
integrals),
\item $TU=\mathsf{Ker}\left(\left.\Pi\right|_{U}\right)_{*}\oplus\mathsf{Ker}F_{*}$
(transversality),
\item $\mathsf{Ker}F_{*}\subseteq\left(\mathsf{Ker}F_{*}\right)^{\bot}$
(isotropy).
\end{itemize}
In other words, from $\Sigma$ we have, around every point of $M$,
a first integrals submersion, transverse to $\Pi$ and isotropic {[}see
Eqs. \eqref{fi}, \eqref{tr} and \eqref{is}{]}. $U$ and $F$ are
given by the formulae
\begin{equation}
U\coloneqq\Sigma\left(V\right)\;\;\;\textrm{and \;\;\;}F\coloneqq p_{\Lambda}\circ\left(\left.\Sigma\right|_{V}\right)^{-1},\label{cf}
\end{equation}
where $V\subseteq N\times\Lambda$ is an open subset for which $\left.\Sigma\right|_{V}$
is a diffeomorphism with its image and $p_{\Lambda}:N\times\Lambda\rightarrow\Lambda$
is the canonical projection onto the second factor.

Reciprocally (see also \cite{gp}), from a submersion $F:M\rightarrow\Lambda$
satisfying \eqref{fi}, \eqref{tr} and \eqref{is}, we can construct,
around every point of $M$, a neighborhood $U$ and a local isotropic
complete solution $\Sigma$ of the $\Pi$-HJE. The involved subset
$U$ is one for which $\left.\left(\Pi,F\right)\right|_{U}$ is a
diffeomorphism with its image, and $\Sigma$ is given by
\begin{equation}
\Sigma=\left[\left.\left(\Pi,F\right)\right|_{U}\right]^{-1}:\Pi\left(U\right)\times F\left(U\right)\rightarrow U.\label{fc}
\end{equation}

\bigskip{}

Summarizing, an isotropic complete solution gives rise to local isotropic
first integrals \textit{via} the Eq. \eqref{cf}, and isotropic first
integrals give rise to a local isotropic complete solution \textit{via}
the Eq. \eqref{fc}.

\section{The existence theorems}

\label{cond}

Let us continue with a symplectic manifold $\left(M,\omega\right)$
of dimension $d=2s$, a function $H:M\rightarrow\mathbb{R}$ and a
fibration $\Pi:M\rightarrow N$ with $\mathsf{dim}N=k\geq1$. Now
let us assume that, for a given $m\in M$,
\begin{description}
\item [{i.}] $X_{H}\left(m\right)\notin\mathsf{Ker}\Pi_{*,m}$ (in particular,
$X_{H}\left(m\right)\neq0$);
\item [{ii.}] $\mathsf{Ker}\Pi_{*,m}$ is co-isotropic, i.e. $\left(\mathsf{Ker}\Pi_{*,m}\right)^{\bot}\subseteq\mathsf{Ker}\Pi_{*,m}$.
\end{description}
Note that, defining $l\coloneqq\dim\mathsf{Ker}\Pi_{*,m}=d-k$, condition
\textbf{$\left(\mathbf{ii}\right)$} implies the inequalities
\[
1\leq k\leq s\leq l\leq d.
\]

The next proposition constitutes the central result of the paper.
From it, we shall derive almost all the other results.

\begin{prop} \label{p1} Under above assumptions, the following assertions
hold:
\begin{enumerate}
\item There exists an open neighborhood $U$ of $m$ and a distribution
$\mathcal{D}$ along $U$ such that:
\begin{enumerate}
\item $\mathsf{Im}\left.X_{H}\right|_{U}\subset\mathcal{D}$,
\item $TU=\mathsf{Ker}\left(\left.\Pi\right|_{U}\right)_{*}\oplus\mathcal{D}$,
\item $\mathcal{D}\subseteq\mathcal{D}^{\bot}$,
\item $\mathcal{D}$ is integrable. In particular, there exists a submersion
$F:U\rightarrow\mathbb{R}^{l}$ such that $\mathcal{D}=\mathsf{Ker}F_{*}$.
\end{enumerate}
\item $\mathcal{D}$ can be chosen such that $\mathcal{D}^{\bot}$ is also
integrable. \\ In other words, a first integrals submersion, transverse
to $\Pi$, isotropic and symplectically complete, can be defined around
$m$.
\end{enumerate}
\end{prop}

\textit{Proof. }Let us begin with the first assertion. To do that,
given $r\in\mathbb{N}$, consider the property \textbf{P}$_{r}$:
there exists a neighborhood $U_{r}$ of $m$ and Hamiltonian vector
fields $X_{1},...,X_{r}\in\mathfrak{X}\left(U_{r}\right)$ which are
linearly independent, mutually orthogonal w.r.t. $\omega$ (and consequently
$\left[X_{i},X_{j}\right]=0$ for all $i,j$),\footnote{Recall that, given to functions $f,g$, for their Hamiltonian vector
fields {[}using convention \eqref{conv}{]} we have that
\[
\left[X_{f},X_{g}\right]=-X_{\omega\left(X_{f},X_{g}\right)}.
\]
Then, if $\omega\left(X_{f},X_{g}\right)=0$, it follows that $\left[X_{f},X_{g}\right]=0$.} with one of them equal to $\left.X_{H}\right|_{U_{r}}$, and such
that the distribution $\mathcal{D}_{r}\coloneqq\left\langle X_{1},...,X_{r}\right\rangle $
satisfies
\begin{equation}
\left.\mathcal{D}_{r}\right|_{p}\cap\mathsf{Ker}\Pi_{*,p}=\left\{ 0\right\} ,\;\;\;\forall p\in U_{r}.\label{mc}
\end{equation}
We shall show by induction the property \textbf{P}$_{r}$ for all
$r\leq k$. From the validity of \textbf{P}$_{k}$, it is clear that
there exists a distribution $\mathcal{D}_{k}$ satisfying the points
$\left(a\right)$ to $\left(d\right)$ above. 

\bigskip{}

If $r=1$, we can take 
\[
U_{1}\coloneqq\left\{ p\in M:X_{H}\left(p\right)\notin\mathsf{Ker}\Pi_{*,p}\right\} 
\]
and define $X_{1}\coloneqq\left.X_{H}\right|_{U_{1}}$. Suppose that
$k>1$ and that \textbf{P}$_{r}$ is true for some $r<k$, and let
us show that, as a consequence, \textbf{P}$_{r+1}$ is also true.
(Note that $r+1\leq k$). Fix a subset $U_{r}$ and Hamiltonian vector
fields $X_{1},...,X_{r}\in\mathfrak{X}\left(U_{r}\right)$ ensured
by the property \textbf{P}$_{r}$. Assume for simplicity that $X_{1}=\left.X_{H}\right|_{U_{r}}$.
Since the given vector fields satisfy $\left[X_{i},X_{j}\right]=0$
for all $i,j$ (see the last footnote), there exists a coordinate
chart $\left(U_{r+1},\left(y_{1},...,y_{d}\right)\right)$ such that
$m\in U_{r\text{+1}}\subseteq U_{r}$ and
\[
X_{1}\left(p\right)=\left.\frac{\partial}{\partial y_{1}}\right|_{p},\ldots,X_{r}\left(p\right)=\left.\frac{\partial}{\partial y_{r}}\right|_{p},\;\;\;\forall p\in U_{r+1}.
\]
Assume for simplicity that $U_{r}=U_{r+1}$. Then we can write
\begin{equation}
X_{1}=\frac{\partial}{\partial y_{1}},\ldots,X_{r}=\frac{\partial}{\partial y_{r}}.\label{xdy}
\end{equation}
Note that, for the submersion $F_{r}\coloneqq\left(y_{r+1},...,y_{d}\right):U_{r}\rightarrow\mathbb{R}^{d-r}$,
we have that {[}see the point $\left(d\right)${]}
\[
\mathsf{Ker}\left(F_{r}\right)_{*}=\left\langle \frac{\partial}{\partial y_{1}},...,\frac{\partial}{\partial y_{r}}\right\rangle =\left\langle X_{1},...,X_{r}\right\rangle =\mathcal{D}_{r}.
\]

Now, consider the Hamiltonian vector fields $X_{y_{1}},...,X_{y_{d}}\in\mathfrak{X}\left(U_{r+1}\right)$.
Using the orthogonality property $\omega\left(X_{i},X_{j}\right)=0$,
for all $1\leq i,j\leq r$, and the fact that
\begin{equation}
\omega\left(X_{y_{a}},X_{j}\right)=\left\langle dy_{a},\frac{\partial}{\partial y_{j}}\right\rangle =\delta_{aj},\;\;\;1\leq a\leq d,\;\;\;1\leq j\leq r,\label{wxy}
\end{equation}
if we write $X_{i}$ as a linear combination of the $X_{y_{a}}$'s,
i.e. $X_{i}=\sum_{a=1}^{d}c_{i}^{a}\,X_{y_{a}}$, we have that $c_{i}^{j}=0$
for all $1\leq i,j\leq r$. That is to say,
\begin{equation}
X_{i}=\sum_{a=r+1}^{d}c_{i}^{a}\,X_{y_{a}}.\label{xixy}
\end{equation}
Note that, since $r<k\leq d/2$, we have that $d-r>r$. Then, since
the matrix with coefficients $c_{i}^{a}$'s (with $i=1,...,r$ and
$a=r+1,...,d$) must have maximal rank $r$ (because the $X_{i}$'s
are l.i.), we can reorder the $d-r$ coordinates $y_{a}$'s (with
$a$ between $r+1$ and $d$) to ensure that the last $r$ columns
define an invertible matrix. For this new order, it is easy to see
that the set of vector fields
\begin{equation}
\left\{ X_{1},...,X_{r},X_{y_{1}},...,X_{y_{r}},X_{y_{r+1}},...,X_{y_{d-r}}\right\} \label{basis}
\end{equation}
gives a basis for $TU_{r+1}$. Let us focus on the fields above evaluated
at $m$. We shall see that there exists a natural number $b$ between
$r+1$ and $d-r$ such that
\begin{equation}
\left\{ X_{1}\left(m\right),...,X_{r}\left(m\right),X_{y_{b}}\left(m\right)\right\} \cap\mathsf{Ker}\Pi_{*,m}=\left\{ 0\right\} .\label{eq:ab}
\end{equation}
This would say that the Hamiltonian vector fields $X_{1},...,X_{r},X_{r+1},$
with $X_{r+1}=X_{y_{b}}$, are independent, mutually orthogonal, with
one of them equal to $\left.X_{H}\right|_{U_{r+1}}$, and such that
the distribution generated by them $\mathcal{D}_{r+1}\coloneqq\left\langle X_{1},...,X_{r+1}\right\rangle $
satisfies
\[
\left.\mathcal{D}_{r+1}\right|_{p}\cap\mathsf{Ker}\Pi_{*,p}=\left\{ 0\right\} ,\;\;\;\forall p\in U_{r+1},
\]
shrinking $U_{r+1}$ if necessary {[}since \eqref{eq:ab} is an open
condition{]}. In other words, we would prove that \textbf{P}$_{r+1}$
is true.

\bigskip{}

In order to show the validity of \eqref{eq:ab}, we shall proceed
by reductio ad absurdum. Suppose first that $r<k/2$. If \eqref{eq:ab}
does not hold for any $b$ between $r+1$ and $d-r$, then there exist
numbers $A_{b}^{j},\alpha_{b}\in\mathbb{R}$ such that
\begin{equation}
0\neq u_{b}\coloneqq\sum_{j=1}^{r}A_{b}^{j}\,X_{j}\left(m\right)+\alpha_{b}\,X_{y_{b}}\left(m\right)\in\mathsf{Ker}\Pi_{*,m},\;\;\;\;b=r+1,...,d-r.\label{li}
\end{equation}
If $\sum_{j=1}^{r}A_{b}^{j}\,X_{j}\left(m\right)=0$, then $\alpha_{b}\neq0$
(to ensure that $u_{b}\neq0$), and if $\sum_{j=1}^{r}A_{b}^{j}\,X_{j}\left(m\right)\neq0$,
the inductive hypothesis ensures that $\sum_{j=1}^{r}A_{b}^{j}\,X_{j}\left(m\right)$
is outside $\mathsf{Ker}\Pi_{*,m}$, what forces again to have $\alpha_{b}\neq0$.
Hence, the vectors $u_{b}$'s define a set of
\[
\left(d-r\right)-\left(r+1\right)+1=d-2r
\]
linearly independent vectors inside de $\mathsf{Ker}\Pi_{*,m}$. But
\[
d-2r>d-k=l=\dim\mathsf{Ker}\Pi_{*,m}
\]
(because we are assuming that $2r<k$), which is a contradiction.
Now, let us consider the complementary case: $k/2\leq r<k$. Again,
if \eqref{eq:ab} does not hold for any $b$ between $r+1$ and $d-r$,
we have the $d-2r$ independent vectors $u_{b}\in\mathsf{Ker}\Pi_{*,m}$
as above {[}see \eqref{li}{]}. Since now 
\[
d-2r\leq d-k=l,
\]
to complete a basis for $\mathsf{Ker}\Pi_{*,m}$ we can add $l-\left(d-2r\right)=2r-k$
additional vectors to the $u_{b}$'s. If we call them $\left\{ v_{1},...,v_{2r-k}\right\} $,
the set
\[
\left\{ u_{r+1},...,u_{d-r},v_{1},...,v_{2r-k}\right\} 
\]
is a basis for $\mathsf{Ker}\Pi_{*,m}$. Note that, since $r<k$,
then 
\begin{equation}
2r-k<r.\label{des}
\end{equation}
Let us show that there exists a non-null vector $v\in\left.\mathcal{D}_{r}\right|_{m}$
such that 
\begin{equation}
\omega\left(v,v_{a}\right)=0,\;\;\;1\leq a\leq2r-k.\label{wbb}
\end{equation}
Writing $v=\sum_{i=1}^{r}x_{i}\,X_{i}\left(m\right)$ and, using Eq.
\eqref{basis},
\[
v_{a}=\sum_{i=1}^{r}y_{i}\,X_{i}\left(m\right)+\sum_{j=1}^{d-r}B_{a}^{j}\,X_{y_{j}}\left(m\right),
\]
Eq. \eqref{wbb} translates to {[}see Eq. \eqref{wxy} and recall
that $d-r>r${]}
\[
\mathbf{B}\cdot\left(\begin{array}{c}
x_{1}\\
x_{2}\\
\vdots\\
x_{r}
\end{array}\right)=\mathbf{0}.
\]
Since $\mathbf{B}$ is a rectangular $\left(2r-k\right)\times r$
matrix, with more columns than rows {[}recall Eq. \eqref{des}{]},
the above linear system has a non-trivial solution, as we wanted to
show. Since in addition {[}see \eqref{wxy} and \eqref{li}{]}
\[
\omega\left(v,u_{a}\right)=0,\;\;\;r+1\leq a\leq d-r,
\]
it follows that $v$ belongs to $\left(\mathsf{Ker}\Pi_{*,m}\right)^{\bot}$.
But $\mathsf{Ker}\Pi_{*,m}$ is co-isotropic, what implies that $v\in\mathsf{Ker}\Pi_{*,m}$.
So, $v\neq0$ and $v\in\left.\mathcal{D}_{r}\right|_{m}\cap\mathsf{Ker}\Pi_{*,m}$,
which is a contradiction. As a consequence, for any $r<k$, there
always exists $b$ fulfilling Eq. \eqref{eq:ab}, as claimed.

\bigskip{}

Summarizing, if $r+1=k$, calculations above show that the neighborhood
$U\coloneqq U_{k}$ of $m$ and the distribution $\mathcal{D}\coloneqq\left\langle X_{1},...,X_{k}\right\rangle $
along $U$ satisfy the points $a$, $b$, $c$ and $d$. It rests
to show the last assertion. To do that, fix a coordinate chart $\left(U,\left(y_{1},...,y_{d}\right)\right)$
(shrinking $U$ if needed) for which \eqref{xdy} holds (for $r=k$).
Note that, according to \eqref{wxy},
\[
X_{1},...,X_{k},X_{y_{k+1}},...,X_{y_{d-k}}\in\mathcal{D}^{\bot},
\]
and according to Eq. \eqref{xixy} and the discussion below, we can
reorder the coordinates such that above $d-k$ vector fields are linearly
independent. Finally, since $\mathsf{dim}\mathcal{D}^{\bot}=d-k$,
it follows that
\[
\left\langle X_{1},...,X_{k},X_{y_{k+1}},...,X_{y_{d-k}}\right\rangle =\mathcal{D}^{\bot}.
\]
So, we just need to show for the vector fields $X_{y_{a}},X_{y_{b}}\in\mathcal{D}^{\bot}$
that $\left[X_{y_{a}},X_{y_{b}}\right]\in\mathcal{D}^{\bot}$ (because
for the other pairs of vector fields this follows from orthogonality).
On the one hand, we have that
\[
\omega\left(\left[X_{y_{a}},X_{y_{b}}\right],X_{i}\right)=-\left\{ \left\{ y_{a},y_{b}\right\} ,f_{i}\right\} =\left\{ \left\{ f_{i},y_{a}\right\} ,y_{b}\right\} +\left\{ \left\{ y_{b},f_{i}\right\} ,y_{a}\right\} ,
\]
being $f_{i}$ the Hamiltonian of $X_{i}$ and $\left\{ \cdot,\cdot\right\} $
the Poisson bracket associated to $\omega$. On the other hand, using
\eqref{wxy}, $\left\{ f_{i},y_{a}\right\} =\omega\left(X_{i},X_{y_{a}}\right)=0$
(for $1\leq i\leq k$ and $k+1\leq a\leq d-k$). Consequently $\left[X_{y_{a}},X_{y_{b}}\right]\in\mathcal{D}^{\bot}$,
which implies that $\mathcal{D}^{\bot}$ is integrable. $\;\;\;\triangle$

\subsection{Isotropic complete solutions of the $\Pi$-HJE}

Using above proposition and the duality between first integrals and
complete solutions, the main result of this paper immediately follows.
\begin{thm}
Consider a symplectic manifold $\left(M,\omega\right)$, a Hamiltonian
vector field $X_{H}\in\mathfrak{X}\left(M\right)$ and a fibration
$\Pi:M\rightarrow N$. Assume that $X_{H}\left(m\right)\notin\mathsf{Ker}\Pi_{*,m}$
and $\mathsf{Ker}\Pi_{*,m}$ is co-isotropic for some point $m\in M$.
Then, there exists a neighborhood $U$ of $m$ and a local isotropic
complete solution of the $\Pi$-HJE for $X_{H}$ along $U$.
\end{thm}
\textit{Proof}. From proposition above we know that there exists a
neighborhood $U$ of $m$ and a submersion $F:U\rightarrow\mathbb{R}^{l}$
satisfying \eqref{fi}, \eqref{tr} (replacing $M$ by $U$ and $\Pi$
by $\left.\Pi\right|_{U}$) and \eqref{is}. It is clear from the
transversality property \eqref{tr} that $U$ can be taken such that
$\left(\left.\Pi\right|_{U},F\right)$ is a diffeomorphism. Then,
as we said in Section \ref{dual} {[}see Eq. \eqref{fc}{]}, the inverse
of $\left(\left.\Pi\right|_{U},F\right):U\rightarrow\Pi\left(U\right)\times F\left(U\right)$
is a local isotropic complete solution of the $\Pi$-HJE for $X_{H}$
along $U$.$\;\;\;\triangle$
\begin{rem*}
Note that the second point of Proposition \ref{p1} is not used in
above theorem. It will be used in Section \ref{nci}.

In the following, we shall present some other consequences of Proposition
\eqref{p1}. For instance, if no fibration is fixed beforehand, we
shall see in Section \ref{ff} that an isotropic complete solution
exists around $m$ for some fibration $\pi$, with the only condition
that $X_{H}\left(m\right)\neq0$ (i.e. $m$ is not a critical point
for $X_{H}$).
\end{rem*}

\subsection{Standard complete solutions}

Suppose that we are in the standard situation, i.e. $M=T^{*}Q$, $\omega=\omega_{Q}$,
$N=Q$ and $\Pi=\pi_{Q}$. Observe that, following above notation,
we have that $\mathsf{dim}Q=s=k$. It is well-known that $\mathsf{Ker}\left(\pi_{Q}\right)_{*,m}\subseteq T_{m}T^{*}Q$
is Lagrangian for all $m\in T^{*}Q$ (and consequently co-isotropic)
with respect to $\omega_{Q}$. In fact, fixing a Darboux coordinates
chart $\left(U,\varphi=\left(q^{1},...,q^{s},p_{1},...,p_{s}\right)\right)$
around $m$, since
\[
\mathsf{Ker}\left(\pi_{Q}\right)_{*,m}=\left\langle \left.\frac{\partial}{\partial p_{1}}\right|_{m},...,\left.\frac{\partial}{\partial p_{s}}\right|_{m}\right\rangle ,
\]
the identity 
\[
\mathsf{Ker}\left(\pi_{Q}\right)_{*}=\left(\mathsf{Ker}\left(\pi_{Q}\right)_{*}\right)^{\bot}
\]
is immediate. Then, if $X_{H}\left(m\right)\notin\mathsf{Ker}\left(\pi_{Q}\right)_{*,m}$,
above theorem ensures that a local isotropic complete solution exists
around $m$. But, since $k=s$, such a solution is actually Lagrangian.
As explained in the Introduction (see \cite{gp} for more details),
the Lagrangian complete solutions of the $\pi_{Q}$-HJE for $H$ are
locally given by the expression
\[
\Sigma\left(q,\lambda\right)=\mathsf{d}W_{\lambda}\left(q\right),
\]
with each function $W_{\lambda}$ satisfying \eqref{ehje}. All that
can be condensed in the next result.
\begin{cor}
Given an $s$-manifold $Q$ and a Hamiltonian function $H:T^{*}Q\rightarrow\mathbb{R}$,
around every point $m\in T^{*}Q$ such that $X_{H}\left(m\right)\notin\mathsf{Ker}\left(\pi_{Q}\right)_{*,m}$
there exists a neighborhood $U\subseteq T^{*}Q$ of $m$, another
$s$-manifold $\Lambda$ and a family of functions $W_{\lambda}:\pi_{Q}\left(U\right)\rightarrow\mathbb{R}$
such that 
\[
\mathsf{d}\left(H\circ\mathsf{d}W_{\lambda}\right)=0,\;\;\;\forall\lambda\in\Lambda,
\]
and
\[
\left(q,\lambda\right)\in\pi_{Q}\left(U\right)\times\Lambda\longmapsto\mathsf{d}W_{\lambda}\left(q\right)\in U
\]
is a diffeomorphism. In other terms, under above conditions, a complete
family of Hamilton's characteristic functions exists around $m$.
\end{cor}
Just to analyze a particular situation, but a very common one, suppose
that $H$ is of the form
\[
H=\frac{1}{2}\,\mathfrak{H}+h\circ\pi_{Q},
\]
where $\mathfrak{H}:T^{*}Q\rightarrow\mathbb{R}$ is a quadratic form
defined by a Riemannian metric $\phi:TQ\times_{Q}TQ\rightarrow\mathbb{R}$
on $Q$, i.e.
\[
\mathfrak{H}\left(m\right)=\phi\left(\phi^{\sharp}\left(m\right),\phi^{\sharp}\left(m\right)\right),
\]
 and $h\in C^{\infty}\left(Q\right)$. In other words, suppose that
$H$ is \textit{simple}. Then, as it is easy to show,
\[
\left(\pi_{Q}\right)_{*}\left(X_{H}\left(m\right)\right)=\phi^{\sharp}\left(m\right),
\]
what implies that $X_{H}\left(m\right)\in\mathsf{Ker}\left(\pi_{Q}\right)_{*,m}$
if and only if $m=0$ (i.e. $m$ belongs to the null subbundle of
$T^{*}Q$). As a consequence, we have the following corollary.
\begin{cor}
Consider a manifold $Q$ and a simple Hamiltonian function $H:T^{*}Q\rightarrow\mathbb{R}$.
Then, if $m\neq0$, a complete family of Hamilton's characteristic
functions exists around $m$.
\end{cor}

\subsection{Non-commutative integrability}

\label{nci}

This section is included just for completeness. Let us fix again a
symplectic manifold $\left(M,\omega\right)$ with $\mathsf{dim}M=2s=d$. 
\begin{defn*}
Consider a function $H\in C^{\infty}\left(M\right)$ and a submersion
$F:M\rightarrow\Lambda$. We shall say that the pair $\left(H,F\right)$
is a \textbf{non-commutative }or \textbf{Mischenko-Fomenko integrable
system} \cite{mf} if $F$ is a first integrals submersion, isotropic
and symplectically complete. And we shall say that $\left(H,F\right)$
is a \textbf{commutative }or \textbf{Arnold-Liouville integrable system
}\cite{ar} if in addition $F$ is Lagrangian.
\end{defn*}
As it is well-known, if $\left(H,F\right)$ is a non-commutative integrable
system, then the Hamiltonian system defined by $H$ can be integrated
up to quadratures.
\begin{rem*}
Usually, additional conditions are asked to $F$, as the compactness
and connectedness of their level submanifolds. In such a case, beside
the integrability up to quadratures, we can ensure (without any calculation)
that the trajectories of the system will be given by quasi-periodic
orbits (see \cite{j} for a review on the subject). But we shall not
consider in this paper these additional conditions.
\end{rem*}
We shall show in this section that, around every non-critical point
of $H$, a local submersion $F:U\rightarrow\mathbb{R}^{l}$ always
exists such that $\left(\left.H\right|_{U},F\right)$ is non-commutative
integrable. We could do it by following a similar procedure to that
applied in the proof of the Proposition \ref{p1}. Since now no fibration
is involved, we can consider properties like \textbf{P}$_{r}$, but
where condition \eqref{mc} is not asked. Also, as in the Appendix
of Ref. \cite{g}, to prove that result we can use the Carathéodory-Jacobi-Lie
theorem \cite{lm} (see also the extension given in \cite{eva}).
Nevertheless, to take advantage of Proposition \eqref{p1}, in this
paper we choose to prove it as a corollary of such a proposition,
with the aid of the result below.
\begin{prop}
\label{p2} Given a (non-necessarily Hamiltonian) vector field $X$,
a point $m$ such that $X\left(m\right)\neq0$, and a natural number
$k\leq s$, there exist an open neighborhood $U$ of $m$ and a submersion
$\pi:U\rightarrow\mathbb{R}^{k}$ such that $X\left(m\right)\notin\mathsf{Ker}\pi_{*,m}$
and $\left(\mathsf{Ker}\pi_{*,m}\right)^{\bot}\subseteq\mathsf{Ker}\pi_{*,m}$.
\end{prop}
\textit{Proof}. Let $\left(U,\left(q^{1},...,q^{s},p_{1},...,p_{s}\right)\right)$
be a Darboux coordinate chart around $m$, and let us write $X$$\left(m\right)$
as a linear combination of the related coordinate vector fields, i.e.
\[
X\left(m\right)=\sum_{i=1}^{s}a^{i}\,\left.\frac{\partial}{\partial q^{i}}\right|_{m}+\sum_{j=1}^{s}b_{j}\,\left.\frac{\partial}{\partial p_{j}}\right|_{m}.
\]
Since $X\left(m\right)\neq0$, then 
\begin{enumerate}
\item $a^{i}\neq0$ for some $i$ or 
\item $b_{j}\neq0$ for some $j$. 
\end{enumerate}
It is clear that the second case reduces to the first one if we make
the canonical transformation 
\[
\left(q^{1},...,q^{s},p_{1},...,p_{s}\right)\mapsto\left(-p_{1},...,-p_{s},q^{1},...,q^{s}\right).
\]
So, we can assume that we are in the first case. In such a case, reordering
the coordinates $q^{i}$'s, if needed (and making the same reordering
for the coordinates $p_{i}$'s), we can always assume that $a^{1}\neq0$.
This means that
\[
\left\{ X\left(m\right),\left.\frac{\partial}{\partial q^{2}}\right|_{m},...,\left.\frac{\partial}{\partial q^{s}}\right|_{m},\left.\frac{\partial}{\partial p_{1}}\right|_{m},...,\left.\frac{\partial}{\partial p_{s}}\right|_{m}\right\} 
\]
is a basis for $T_{m}M$. Consider now the submersion $\pi\coloneqq\left(q_{1},...,q_{k}\right):U\rightarrow\mathbb{R}^{k}$.
Since 
\[
\mathsf{Ker}\pi_{*}=\left\langle \frac{\partial}{\partial q^{k+1}},...,\frac{\partial}{\partial q^{s}},\frac{\partial}{\partial p_{1}},...,\frac{\partial}{\partial p_{s}}\right\rangle 
\]
and
\[
\left(\mathsf{Ker}\pi_{*}\right)^{\bot}=\left\langle \frac{\partial}{\partial p_{1}},...,\frac{\partial}{\partial p_{k}}\right\rangle ,
\]
it is clear that $X\left(m\right)\notin\mathsf{Ker}\pi_{*,m}$ and
$\mathsf{Ker}\pi_{*}$ is co-isotropic. Thus, $\pi$ is the submersion
we are looking for.$\;\;\;\triangle$

\bigskip{}

Combining Propositions \ref{p1} and \ref{p2}, we have the wanted
result.
\begin{thm}
\label{ncie} Consider a $2s$-dimensional symplectic manifold $\left(M,\omega\right)$
and a function $H\in C^{\infty}\left(M\right)$. Assume that, for
a given point $m\in M$, $\mathsf{d}H\left(m\right)\neq0$. Then,
for every $l$ such that $2s>l\geq s$, there exists a neighborhood
$U$ of $m$ and a submersion $F:U\rightarrow\mathbb{R}^{l}$ such
that the pair $\left(\left.H\right|_{U},F\right)$ is a non-commutative
integrable system. In particular, choosing $l=s$, we have that $\left(\left.H\right|_{U},F\right)$
is a commutative integrable system.
\end{thm}
\textit{Proof}. If $\mathsf{d}H\left(m\right)\neq0$, then $X_{H}\left(m\right)\neq0$,
and consequently, using Proposition \ref{p2}, for each $k\leq s$
there exists a fibration $\pi:U\rightarrow\pi\left(U\right)\subseteq\mathbb{R}^{k}$
such that $X_{H}\left(m\right)\notin\mathsf{Ker}\pi_{*,m}$ and $\mathsf{Ker}\pi_{*,m}$
is co-isotropic. In this situation, given $l$ such that $2s>l\geq s$,
and taking $k=2s-l$, Proposition \ref{p1} ensures the existence
of a submersion $F:U\rightarrow\mathbb{R}^{l}$ (shrinking $U$ if
needed) satisfying precisely the conditions of the non-commutative
integrability (see $1.\left(d\right)$ and then $1.\left(a\right)$,
$1.\left(c\right)$ and $2$).$\;\;\;\triangle$

\subsection{Isotropic complete solutions for some fibration}

\label{ff}

If no fibration is fixed beforehand, we have the next result. 
\begin{thm}
\label{ffe} Consider a $2s$-dimensional symplectic manifold $\left(M,\omega\right)$
and a function $H\in C^{\infty}\left(M\right)$. Assume that, for
a given point $m\in M$, $\mathsf{d}H\left(m\right)\neq0$. Then,
for every $k\leq s$, there exists a neighborhood $U$ of $m$, a
fibration $\pi:U\rightarrow\pi\left(U\right)\subseteq\mathbb{R}^{k}$
and an isotropic complete solution of the $\pi$-HJE for $\left.X_{H}\right|_{U}$. 
\end{thm}
\textit{Proof}. According to Proposition \ref{p2}, there exist an
open neighborhood $U$ and fibration a $\pi:U\rightarrow\pi\left(U\right)\subseteq\mathbb{R}^{k}$
such that $X_{H}\left(m\right)\notin\mathsf{Ker}\pi_{*,m}$ and $\left(\mathsf{Ker}\pi_{*,m}\right)^{\bot}\subseteq\mathsf{Ker}\pi_{*,m}$.
On the other hand, Proposition \ref{p1} ensures the existence of
another fibration $F:U\rightarrow F\left(U\right)\subseteq\mathbb{R}^{l}$
(shrinking $U$ if needed), with $l=d-k$, such that:
\begin{itemize}
\item $\mathsf{Im}\left.X_{H}\right|_{U}\subset\mathsf{Ker}F_{*}$,
\item $TU=\mathsf{Ker}\pi_{*}\oplus\mathsf{Ker}F_{*}$,
\item $\mathsf{Ker}F_{*}\subseteq\left(\mathsf{Ker}F_{*}\right)^{\bot}$.
\end{itemize}
Then, by duality, assuming for simplicity that $\left(\pi,F\right):U\rightarrow\pi\left(U\right)\times F\left(U\right)$
is a diffeomorphism (otherwise it is enough to change $U$ by a smaller
open set), we know that $\Sigma\coloneqq\left(\pi,F\right)^{-1}$
is an isotropic complete solution of the $\pi$-HJE for $\left.X_{H}\right|_{U}$.
$\;\;\;\triangle$

\bigskip{}
Concluding, around every non-critical point of a Hamiltonian system
there exists:
\begin{itemize}
\item a submersion $F:U\rightarrow\mathbb{R}^{l}$ such that $\left(\left.H\right|_{U},F\right)$
is non-commutative integrable (see Theorem \ref{ncie});
\item an isotropic complete solution $\Sigma$ of the $\pi$-HJE for some
fibration $\pi:U\rightarrow\pi\left(U\right)$ (see Theorem \ref{ffe}). 
\end{itemize}
As we said at the beginning of the paper, above existence results
do not mean that, around a non-critical point, every Hamiltonian system
is exactly solvable. The point is that, in order to ensure exact solvability,
it is not enough to know that the mentioned objects exist, but we
need to have an explicit expression of them.

\section*{Acknowledgements}

The author thanks CONICET for its financial support.

\end{document}